\RequirePackage{ifpdf}
\ifpdf 
\documentclass[pdftex]{sigma}
\else
\documentclass{sigma}
\fi

\begin{document}
\allowdisplaybreaks

\renewcommand{\thefootnote}{$\star$}

\renewcommand{\PaperNumber}{022}

\FirstPageHeading

\ShortArticleName{Conformal Killing--Yano Tensors on Manifolds with
Mixed 3-Structures}

\ArticleName{Conformal Killing--Yano Tensors on Manifolds \\ with
Mixed 3-Structures\footnote{This paper is a contribution to the Proceedings of the XVIIth International Colloquium on Integrable Systems and Quantum Symmetries (June 19--22, 2008, Prague, Czech Republic). The full collection
is available at
\href{http://www.emis.de/journals/SIGMA/ISQS2008.html}{http://www.emis.de/journals/SIGMA/ISQS2008.html}}}

\Author{Stere IANU\c{S}~$^\dag$, Mihai VISINESCU~$^\ddag$ and
Gabriel Eduard V\^{I}LCU~$^{\dag\S}$}

\AuthorNameForHeading{S. Ianu\c{s}, M. Visinescu and G.E. V\^{\i}lcu}

\Address{$^\dag$~University of Bucharest, Faculty of Mathematics and
Computer Science,\\
\hphantom{$^\dag$}~Str.\ Academiei, Nr.~14, Sector~1, Bucharest 70109,
Romania} \EmailD{\href{mailto:ianus@gta.math.unibuc.ro}{ianus@gta.math.unibuc.ro}}

\Address{$^\ddag$~National Institute for Physics and Nuclear
        Engineering, Department of Theoretical Physics,\\
\hphantom{$^\ddag$}~P.O.~Box M.G.-6,         Magurele, Bucharest, Romania}
\EmailD{\href{mailto:mvisin@theory.nipne.ro}{mvisin@theory.nipne.ro}}
\URLaddressD{\url{http://www.theory.nipne.ro/~mvisin/}}

\Address{$^\S$~Petroleum-Gas University of Ploie\c sti,
         Department of Mathematics and Computer Science,\\
\hphantom{$^\S$}~Bulevardul Bucure\c sti, Nr.~39, Ploie\c sti 100680, Romania}
\EmailD{\href{mailto:gvilcu@mail.upg-ploiesti.ro}{gvilcu@mail.upg-ploiesti.ro}}

\ArticleDates{Received October 30, 2008, in f\/inal form February 16,
2009; Published online February 23, 2009}

\Abstract{We show the existence of conformal Killing--Yano tensors on
a manifold endowed with a mixed 3-Sasakian structure.}


\Classification{53C15; 81T20}

\section{Introduction}

Investigations of the properties of space-times of higher dimensions
($D > 4$) have recently attracted considerable attention as a result
of their appearance in theories of unif\/ication such as string and
$M$ theories.

It is important to study the symmetries of these space-times in
order to characterize the geodesic motions in them. In general a
space-time could possess explicit and {\it hidden} symmetries
encoded in the multitude of Killing vectors and higher Killing
tensors respectively. The symmet\-ries allow one to def\/ine conserved
quantities along geodesics linear and polynomial in canonical
momenta. Their existence guarantees the integrability of the
geodesic motions and is intimately related to separability of
Hamilton--Jacobi (see, e.g.~\cite{BEN}) and the Klein--Gordon
equation \cite{CARTER} at the quantum level.

The next most simple objects that can be studied in connection with
the symmetries of a manifold after Killing tensors are Killing--Yano
tensor~\cite{YAN}. Their physical interpretation remained obscure
until Floyd \cite{FLO} and Penrose \cite{PEN} showed that the
Killing tensor $K^{\mu\nu}$ of the $4$-dimensional Kerr--Newman
space-time admits a certain square-root which def\/ines a Killing--Yano
tensor. Subsequently it was realized \cite{GRH} that a Killing--Yano
tensor generate additional {\it supercharges} in the dynamics of
pseudo-classical spinning particles. In this way it was realized the
natural connection between Killing--Yano tensors and {\it
supersymmetries} \cite{PLY}. Passing to quantum Dirac equation it
was discovered \cite{CARMCL} that Killing--Yano tensors generate
conserved {\it non-standard} Dirac operators which commute with the
{\it standard} one.

The conformal extension of the Killing tensor equation determines
the conformal Killing tensors \cite{JJ} which def\/ine f\/irst integrals
of the null geodesic equation. Investigation of the properties of
higher-dimensional space-times has pointed out the role of the
conformal Killing--Yano tensors to generate background metrics with
black-hole solutions (see, e.g.~\cite{FRO}).

The aim of this paper is to investigate the existence of conformal
Killing--Yano tensors on some manifolds endowed with special
structures which could be relevant in the theories of modern physics
\cite{CAR,GRH,LP,TAN,VV}.

Versions of $M$-theory could be formulated in a space-time with
various number of time dimensions giving rise to exotic space-time
signatures. The $M$-theory in $10+1$ dimensions is linked via
dualities to a $M^{*}$ theory in $9+2$ dimensions and a $M^{'}$
theory in $6+5$ dimensions. Various  limits of these will give rise
to $IIA$- and $IIB$-like string theories in many variants of
dimensions and signatures~\cite{HUL}.

The paraquaternionic structures arise in a natural way in
theoretical physics, both in string theory and integrable systems
\cite{CMMS,DW,H,OV}. The counterpart in odd dimension of a
paraquaternionic structure was introduced in~\cite{IMV}. It is
called mixed 3-structure, which appears in a natural way on
lightlike hypersurfaces in paraquaternionic manifolds. A compatible
metric with a mixed 3-structures is necessarily semi-Riemann and
mixed 3-Sasakian manifolds are Einstein \cite{CP,IV}, hence the
possible importance of these structures in theoretical physics.

\section[Killing vector fields and their generalizations]{Killing vector f\/ields and their generalizations}

Let $(M,g)$ be a semi-Riemannian manifold. A vector f\/ield $X$ on $M$
is said to be a Killing vector f\/ield if the Lie derivative with
respect to $X$ of the metric $g$ vanishes:
\begin{gather}\label{1}
L_Xg=0 .
\end{gather}
Alternatively, if $\nabla$ denotes the Levi-Civita connection of
$g$, then (\ref{1}) can be rewritten as
\[
g(\nabla_YX,Z)+g(Y,\nabla_ZX)=0 ,
\]
for all vector f\/ields $Y$, $Z$ on $M$.

It is easy to see that the Lie bracket of two Killing f\/ields is
still a Killing f\/ield and thus the Killing f\/ields on a manifold $M$
form a Lie subalgebra of vector f\/ields on $M$. This is the Lie
algebra of the isometry group of the manifold $M$, provided that $M$
is compact.

Killing vector f\/ields can be generalized to conformal Killing vector
f\/ields \cite{YAN}, i.e.\ vector f\/ields with a f\/low preserving a given
conformal class of metrics. A natural generalization of conformal
Killing vector f\/ields is given by the conformal Killing--Yano tensors
\cite{KSW}. Roughly speaking, a~conformal Killing--Yano tensor of
rank $p$ on a semi-Riemannian manifold $(M,g)$ is a $p$-form
$\omega$ which satisf\/ies:
\begin{gather}\label{4}
\nabla_X\omega=\frac{1}{p+1}X \lrcorner
d\omega-\frac{1}{n-p+1}X^*\wedge d^*\omega  ,
\end{gather}
for any vector f\/ield $X$ on $M$, where $\nabla$ is the Levi-Civita
connection of $g$, $n$ is the dimension of~$M$,~$X^*$ is the 1-form
dual to the vector f\/ield $X$ with respect to the metric~$g$,
$\lrcorner$ is the operator dual to the wedge product and~$d^*$ is
the adjoint of the exterior derivative $d$. If $\omega$ is co-closed
in~(\ref{4}), then we obtain the def\/inition of a Killing--Yano tensor
(introduced by Yano~\cite{YAN}). We can easily see that for $p=1$,
they are dual to Killing vector f\/ields.

We remark that Killing--Yano tensors are also called Yano tensors or
Killing forms, and conformal Killing--Yano tensors are sometimes
referred as conformal Yano tensors, conformal Killing forms or
twistor forms \cite{BMS, MS, SEM}.

For generalizations of the Killing vectors one might also consider
higher order symmetric tensors. Let $\rho$ be a covariant symmetric
tensor f\/ield of rank $r$ on a semi-Riemannian manifold $(M,g)$. Then
$\rho$ is a Killing tensor f\/ield if the symmetrization of covariant
derivative of $\rho$ vanishes identically:
\[
\nabla_{(\lambda} \rho_{\mu_1\dots \mu_r)} = 0 .
\]
If the Ricci tensor of a semi-Riemannian manifold $(M,g)$ is
Killing, then $(M,g)$ is an Einstein-like space~\cite{GR}.

The relevance in physics of the Killing tensors is given by the
following proposition which could be easily proved:
\begin{proposition}
A symmetric tensor $\rho$ on $M$ is a Killing tensor if and only if
the quantity
\[
K = \rho_{\mu_1\dots \mu_r} \dot{c}^{\mu_1}\cdots \dot{c}^{\mu_r} 
\]
is constant along every geodesic $c$ in $M$.
\end{proposition}

Here the over-dot denotes the ordinary proper time derivative and
the proposition ensures that $K$ is a f\/irst integral of the geodesic
equation.

These two generalizations  of the Killing vector equation could be
related. Let $\omega_{\mu_1\dots \mu_p}$ be a~Killing--Yano tensor, then
the tensor f\/ield
\[
\rho_{\mu\nu} =  \omega_{\mu\mu_2\dots \mu_p}
\omega_{\;\;\;\;\;\;\;\;\;\;\nu}^{\mu_2\dots \mu_p} 
\]
is a Killing tensor and it sometimes refers to it as the associated
tensor with $\omega$. However, the converse statement is not true in
general: not all Killing tensors of rank $2$ are associated with a
Killing--Yano tensor. That is the case of the spaces which admit
Killing tensors but no Killing--Yano tensors.

\section{Manifolds with mixed 3-structures}

Let $M$ be a dif\/ferentiable manifold equipped with a triple
$(\phi,\xi,\eta)$, where $\phi$ is a a f\/ield  of endomorphisms of
the tangent spaces, $\xi$ is a vector f\/ield and $\eta$ is a 1-form
on $M$ such that:
\[
\phi^2=-\epsilon I+\eta\otimes\xi,\qquad \eta(\xi)=\epsilon .
\]
If $\epsilon=1$ then $(\phi,\xi,\eta)$ is said to be an almost
contact structure on $M$ (see \cite{BLR}), and if $\epsilon=-1$ then
$(\phi,\xi,\eta)$ is said to be an almost paracontact structure on
$M$ (see~\cite{MTS}).

\begin{definition}[\cite{IMV}]
Let $M$ be a dif\/ferentiable manifold which admits an almost contact structure
$(\phi_1,\xi_1,\eta_1)$ and two almost paracontact
structures $(\phi_2,\xi_2,\eta_2)$ and
$(\phi_3,\xi_3,\eta_3)$, satisfying the following
conditions:
\begin{gather*}
\eta_\alpha(\xi_\beta)=0, \qquad \forall \, \alpha\neq\beta ,\\
\phi_\alpha(\xi_\beta)=-\phi_\beta(\xi_\alpha)=\epsilon_\gamma\xi_\gamma ,\\
\eta_\alpha\circ\phi_\beta=-\eta_\beta\circ\phi_\alpha=
\epsilon_\gamma\eta_\gamma,\\
\phi_\alpha\phi_\beta-\eta_\beta\otimes\xi_\alpha=
-\phi_\beta\phi_\alpha+\eta_\alpha\otimes\xi_\beta=
\epsilon_\gamma\phi_\gamma ,
\end{gather*}
where $(\alpha,\beta,\gamma)$ is an even
permutation of (1,2,3) and $\epsilon_1=1$, $\epsilon_2=\epsilon_3=-1$ .

Then the manifold $M$ is said to have a mixed 3-structure
$(\phi_\alpha,\xi_\alpha,\eta_\alpha)_{\alpha=\overline{1,3}}$.
\end{definition}

\begin{definition} If a manifold $M$ with a mixed
3-structure
$(\phi_\alpha,\xi_\alpha,\eta_\alpha)_{\alpha=\overline{1,3}}$
 admits a semi-Riemannian metric $g$ such that:
\begin{gather}\label{11}
g(\phi_\alpha X, \phi_\alpha Y)=\epsilon_\alpha
g(X,Y)-\eta_\alpha(X)\eta_\alpha(Y),\qquad g(X,\xi_\alpha)=\eta_\alpha(X) ,
\end{gather}
for all $X,Y\in\Gamma(TM)$ and $\alpha=1,2,3$, then we say that $M$
has a metric mixed 3-structure and $g$ is called a compatible
metric. Moreover, if $(\phi_1,\xi_1,\eta_1,g)$ is a Sasakian
structure, i.e.\ (see~\cite{BLR}):
\begin{gather}\label{12}
(\nabla_X\phi_1) Y=g(X,Y)\xi_1-\eta_1(Y)X ,
\end{gather}
and $(\phi_2,\xi_2,\eta_2,g)$, $(\phi_3,\xi_3,\eta_3,g)$ are
LP-Sasakian structures, i.e.\ (see~\cite{MTS}):
\begin{gather}\label{13}
(\nabla_X\phi_2) Y=g(\phi_2X,\phi_2Y)\xi_2+\eta_2(Y)\phi_2^2X ,
\\
\label{14}
(\nabla_X\phi_3) Y=g(\phi_3X,\phi_3Y)\xi_3+\eta_3(Y)\phi_3^2X ,
\end{gather}
then
$((\phi_\alpha,\xi_\alpha,\eta_\alpha)_{\alpha=\overline{1,3}},g)$
is said to be a mixed Sasakian 3-structure on $M$.
\end{definition}

We remark that if
$(M,(\phi_\alpha,\xi_\alpha,\eta_\alpha)_{\alpha=\overline{1,3}},g)$
is a manifold with a metric mixed 3-structure then the signature of
$g$ is $(2n+1,2n+2)$ and the dimension of the manifold $M$ is $4n+3$
because one can check that, at each point of $M$, always there
exists a pseudo-orthonormal frame f\/ield given by $\{(E_i,\phi_1 E_i,
\phi_2 E_i, \phi_3 E_i)_{i=\overline{1,n}} , \xi_1, \xi_2,
\xi_3\} .$

The main property of a manifold endowed with a mixed $3$-Sasakian
structure is given by the following theorem (see \cite{CP,IV}):

\begin{theorem} Any $(4n+3)$-dimensional manifold endowed with a
mixed $3$-Sasakian structure is an Einstein space with Einstein
constant $\lambda=4n+2$.
\end{theorem}

Concerning the symmetric Killing tensors let us note that D.E.~Blair
studied in \cite{BL} the almost contact manifold with Killing
structure tensors. He assumed that $M$ has an almost contact metric
structure $(\phi,\xi,\eta,g)$ such that $\phi$ and $\eta$ are
Killing. Then he proved that if $(\phi,\xi,\eta,g)$ is normal, i.e.\
the Nijenhuis tensor $N_\phi$ satisf\/ies (see \cite{BLR}):
\[N_\phi+2d\eta\otimes\xi=0,\] then $(\phi,\xi,\eta,g)$ is a
cosymplectic structure, i.e.\ $\Omega$ and $\eta$ are closed, where
$\Omega$ is the fundamental 2-form of the manifold def\/ined by
\[\Omega(X,Y)=g(\phi X,Y),\] for all vector f\/ields $X,Y$ on $M$.

For a mixed 3-structure with a compatible semi-Riemannian metric
$g$, we have the following result.
\begin{proposition}\label{3.4}
Let $(M,g)$ be a semi-Riemannian manifold. If $(M,g)$ has a mixed
$3$-Sasakian structure
$(\phi_\alpha,\xi_\alpha,\eta_\alpha)_{\alpha=\overline{1,3}}$, then
$(\phi_\alpha)_{\alpha=\overline{1,3}}$ cannot be Killing--Yano
tensor fields.
\end{proposition}
\begin{proof}
If $(\phi_1,\xi_1,\eta_1,g)$ is a Sasakian structure, from
(\ref{12}) we obtain:
\[
(\nabla_X\phi_1)X=g(X,X)\xi_1\neq0\,,
\]
for any non-lightlike vector f\/ield $X$ orthogonal to $\xi_1$.

For LP-Sasakian structures $(\phi_2,\xi_2,\eta_2,g)$ and
$(\phi_3,\xi_3,\eta_3,g)$, from (\ref{13}) and (\ref{14}) we have:
\[
(\nabla_X\phi_\alpha)X=g(\phi_\alpha X,\phi_\alpha X)\xi_\alpha\neq0,
\qquad \alpha\in\{2,3\} ,
\]
for any non-lightlike vector f\/ield $X$ orthogonal to $\xi_2$.
\end{proof}
\begin{theorem}
Let $(M,g)$ be a semi-Riemannian manifold. If $(M,g)$ admits a mixed
$3$-Sasakian structure, then any conformal Killing vector field on
$(M,g)$ is a Killing vector field.
\end{theorem}
\begin{proof} A vector f\/ield $X$ on $M$ is conformal Killing if\/f
\begin{gather}\label{18}
L_Xg=f\cdot g ,
\end{gather}
for $f\in C^\infty(M,\mathbb{R})$.

From (\ref{18}) we have
\[
(L_Xg)(\xi_\alpha,\xi_\alpha)=fg(\xi_\alpha,\xi_\alpha)=\epsilon_\alpha
f .
\]

But, by the Lie operator's properties,
\begin{gather*}
(L_Xg)(\xi_\alpha,\xi_\alpha) = Xg(\xi_\alpha,\xi_\alpha)-
2g(L_X\xi_\alpha,\xi_\alpha)
 = -2g([X,\xi_\alpha],\xi_\alpha)\\
 \phantom{(L_Xg)(\xi_\alpha,\xi_\alpha)}{} = -2g(\nabla_X\xi_\alpha,\xi_\alpha)+2g(\nabla_{\xi_\alpha}X,\xi_\alpha)
  = 2\epsilon_\alpha g(\phi_\alpha X,\xi_\alpha)-
2g(X,\nabla_{\xi_\alpha}\xi_\alpha)
 = 0 ,
\end{gather*}
because $\phi_\alpha X \perp \xi_\alpha$ $(\alpha\in\{1,2,3\})$ and
$\nabla_{\xi_\alpha}\xi_\alpha=0$.

Consequently,
\[
f=\epsilon_\alpha (L_Xg)(\xi_\alpha,\xi_\alpha)=0 ,
\]
so that $L_Xg=0$, i.e.\ $X$ is Killing vector f\/ield. \end{proof}
Actually this result holds on Sasakian, LP-Sasakian and 3-Sasakian
manifolds with the same proof.

\section[Killing vector fields on manifolds with mixed 3-structures]{Killing vector f\/ields on manifolds with mixed 3-structures}

An almost para-hypercomplex structure on a smooth manifold $M$ is a
triple $H=(J_{\alpha})_{\alpha=\overline{1,3}}$, where $J_1$ is an
almost complex structure on $M$ and $J_2$, $J_3$ are almost product
structures on~$M$, satisfying: $J_1J_2J_3=-{\rm Id}$. In this case $(M,H)$
is said to be an almost para-hypercomplex manifold.

A semi-Riemannian metric $g$ on $(M,H)$ is said to be
para-hyperhermitian if it satisf\/ies:
\begin{gather}\label{15}
g(J_\alpha X,J_\alpha Y)=\epsilon_{\alpha}
g(X,Y),\qquad \alpha\in\{1,2,3\} ,
\end{gather}
for all $X,Y\in\Gamma(TM)$, where $\epsilon_1=1$,
$\epsilon_2=\epsilon_3=-1$. In this case, $(M,g,H)$ is called an
almost para-hyperhermitian manifold. Moreover, if each $J_\alpha$ is
parallel with respect to the Levi-Civita connection of $g$, then
$(M,g,H)$ is said to be a para-hyper-K\"{a}hler manifold.

\begin{theorem}\label{4.1} Let $(M,g)$ be a semi-Riemannian manifold.
Then the following five assertions are mutually equivalent:
\begin{enumerate}\itemsep=0pt
\item[$(i)$] $(M,g)$ admits a mixed $3$-Sasakian structure.

\item[$(ii)$] The cone $(C(M),\overline{g})=(M\times\mathbb{R}_+,dr^2+r^2g)$
admits a para-hyper-K\"{a}hler structure.

\item[$(iii)$] There exists three orthogonal Killing vector fields
$\{\xi_1,\xi_2,\xi_3\}$ on $M$, with $\xi_1$ unit spacelike vector
field and $\xi_2$, $\xi_3$ unit timelike vector fields satisfying
\begin{gather}\label{17}
[\xi_\alpha,\xi_\beta]=
(\epsilon_\alpha+\epsilon_\beta)\epsilon_\gamma\xi_\gamma ,
\end{gather}
where $(\alpha,\beta,\gamma)$ is an even permutation of $(1,2,3)$ and
$\epsilon_1=1$, $\epsilon_2=\epsilon_3=-1$, such that the tensor
fields $\phi_\alpha$ of type $(1,1)$, defined by:
\[
\phi_\alpha X=-\epsilon_{\alpha}\nabla_X\xi_\alpha,\qquad
\alpha\in\{1,2,3\} ,
\]
satisfies the conditions \eqref{12}, \eqref{13} and \eqref{14}.

\item[$(iv)$] There exists three orthogonal Killing vector fields
$\{\xi_1,\xi_2,\xi_3\}$ on $M$, with $\xi_1$ unit spacelike vector
field and $\xi_2$, $\xi_3$ unit timelike vector fields satisfying
\eqref{17}, such that:
\[
R(X,\xi_\alpha)Y=g(\xi_\alpha,Y)X-g(X,Y)\xi_\alpha,\qquad
\alpha\in\{1,2,3\} ,
\]
where $R$ is the Riemannian curvature tensor of the Levi-Civita
connection $\nabla$ of $g$.

\item[$(v)$] There exists three orthogonal Killing vector fields
$\{\xi_1,\xi_2,\xi_3\}$ on $M$, with $\xi_1$ unit spacelike vector
field and $\xi_2$, $\xi_3$ unit timelike vector fields satisfying
\eqref{17}, such that the sectional curvature of every section
containing $\xi_1$, $\xi_2$ or $\xi_3$ equals $1$.
\end{enumerate}
\end{theorem}

\begin{proof}
$(i)\Rightarrow(ii)$.
If $(M^{4n+3},(\phi_\alpha,\xi_\alpha,\eta_\alpha)_{\alpha=
\overline{1,3}},g)$ is  a manifold
endowed with a mixed 3-Sasakian structure, then we can def\/ine a
para-hyper-K\"{a}hler structure
$\{J_\alpha\}_{\alpha=\overline{1,3}}$ on the cone
$(C(M),\overline{g})=(M\times\mathbb{R}_+,dr^2+r^2g)$, by:
\begin{gather}\label{18b}
J_\alpha X=\phi_\alpha X-\eta_\alpha(X)\Phi ,\qquad
J_\alpha\Phi=\xi_\alpha ,
\end{gather}
for any $X\in\Gamma(TM)$ and $\alpha\in\{1,2,3\}$, where
$\Phi=r\partial_r$ is the Euler f\/ield on $C(M)$ (see also~\cite{BG}).

$(ii)\Rightarrow(i)$. If the cone
$(C(M),\overline{g})=(M\times\mathbb{R}_+,dr^2+r^2g)$ admits a
para-hyper-K\"{a}hler structure
$\{J_\alpha\}_{\alpha=\overline{1,3}}$, then we can identify $M$
with $M\times\{1\}$ and we have a mixed 3-Sasakian structure
$((\phi_\alpha,\xi_\alpha,\eta_\alpha)_{\alpha=\overline{1,3}},g)$
on $M$ given by:
\begin{gather}\label{19}
\xi_\alpha=J_\alpha(\partial_r),\qquad \phi_\alpha
X=-\epsilon_{\alpha}\nabla_X \xi_\alpha,\qquad
\eta_\alpha(X)=g(\xi_\alpha,X) ,
\end{gather}
for any $X\in\Gamma(TM)$ and $\alpha\in\{1,2,3\}$.

$(ii)\Rightarrow(iii)$. If the cone $C(M)$ admits a
para-hyper-K\"{a}hler structure
$\{J_\alpha\}_{\alpha=\overline{1,3}}$, then using the above
identif\/ications, we can view all vector f\/ields $X$, $Y$ on~$M$ as
vector f\/ields on the $C(M)$. If $\overline{\nabla}$ is the
Levi-Civita connection of $\overline{g}$, then we have the following
formulas (see \cite{ON}):
\begin{gather}\label{20}
\overline{\nabla}_XY=\nabla_XY-rg(X,Y)\partial_r ,\qquad
\overline{\nabla}_{\partial_r}X=\overline{\nabla}_X\partial_r=
\frac{1}{r}X  ,\qquad
\overline{\nabla}_{\partial_r}\partial_r=0 .
\end{gather}

On the other hand, because $(ii)\Rightarrow(i)$, we have that
$((\phi_\alpha,\xi_\alpha,\eta_\alpha)_{\alpha=\overline{1,3}},g)$
given by (\ref{19}) is a~mixed 3-Sasakian structure on $M$ and the
conditions (\ref{12}), (\ref{13}) and (\ref{14}) are satisf\/ied.

Then, using (\ref{15}), (\ref{19}) and (\ref{20}), since $J_1$, $J_2$
and $J_3$ are parallel, we obtain for all vector f\/ields~$X$,~$Y$ on~$M$
and $\alpha\in\{1,2,3\}$:
\begin{gather*}
g(\nabla_X\xi_\alpha,Y)+g(X,\nabla_Y\xi_\alpha)
 = \overline{g}(\overline{\nabla}_X\xi_\alpha+rg(X,\xi_\alpha)\partial_r,Y)
+\overline{g}(X,\overline{\nabla}_Y\xi_\alpha+rg(Y,\xi_\alpha)\partial_r)
\\
\qquad{} = \overline{g}(\overline{\nabla}_X(J_\alpha\partial_r),Y)+
\overline{g}(X,\overline{\nabla}_Y(J_\alpha\partial_r))
 = \overline{g}(J_\alpha\overline{\nabla}_X\partial_r,Y)+
\overline{g}(X,J_\alpha\overline{\nabla}_Y\partial_r)\\
\qquad{} = \frac{1}{r}\left(\overline{g}(J_\alpha
X,Y)+\overline{g}(X,J_\alpha Y)\right)\nonumber = 0 ,
\end{gather*}
and consequently $\xi_\alpha$ is a Killing vector f\/ield, for
$\alpha=1,2,3$.

From (\ref{19}) we also obtain:
\[
g(\xi_\alpha,\xi_\alpha)=\eta_\alpha(\xi_\alpha)=\epsilon_\alpha,\qquad
\alpha\in\{1,2,3\} ,
\]
and
\[
g(\xi_\alpha,\xi_\beta)=\eta_\beta(\xi_\alpha)=0,\qquad
\forall\, \alpha\neq\beta .
\]
Consequently, $\{\xi_1,\xi_2,\xi_3\}$ are
orthogonal Killing vector f\/ields on $M$, with $\xi_1$ unit spacelike
vector f\/ield and $\xi_2$, $\xi_3$ unit timelike vector f\/ields.

Moreover we have:{\samepage
\begin{gather*}
[\xi_\alpha,\xi_\beta] = \nabla_{\xi_\alpha}\xi_\beta-
\nabla_{\xi_\beta}\xi_\alpha
 = -\epsilon_\beta\phi_\beta\xi_\alpha+
\epsilon_\alpha\phi_\alpha\xi_\beta
 = (\epsilon_\alpha+\epsilon_\beta)\epsilon_\gamma\xi_\gamma ,
\end{gather*}
for any even permutation $(\alpha,\beta,\gamma)$ of (1,2,3).}

$(iii)\Rightarrow(ii)$. If $\{\xi_1,\xi_2,\xi_3\}$ are orthogonal
Killing vector f\/ields on $M$, such that $\xi_1$ is unit spacelike,
$\xi_2$, $\xi_3$ are unit timelike vector f\/ields and relations
(\ref{12}), (\ref{13}) and (\ref{14}) are satisf\/ied, then we can
def\/ine a para-hyper-K\"{a}hler structure on the cone
$(C(M),\overline{g})=(M\times\mathbb{R}_+,dr^2+r^2g)$, by~(\ref{18b}), where $\eta_\alpha$ is the 1-form dual to $\xi_\alpha$,
for $\alpha\in\{1,2,3\}$.

$(iii)\Leftrightarrow(iv)$. This equivalence follows from
straightforward computations using the expression of the Riemannian
curvature tensor $R$.

$(iv)\Leftrightarrow(v)$. This equivalence is a simple calculation,
using the formula of the sectional curvature:
\[
K(\pi)=\frac{R(X,Y,X,Y)}{g(X,X)g(Y,Y)-g(X,Y)^2} ,
\]
for a non-degenerate 2-plane $\pi$ spanned by $\{X,Y\}$.
\end{proof}

\begin{corollary}\label{4.2} Let
$(M^{4n+3},(\phi_\alpha,\xi_\alpha,\eta_\alpha)_{\alpha=\overline{1,3}},g)$
be a manifold endowed with a mixed $3$-Sasakian structure. Then:
\begin{enumerate}\itemsep=0pt
\item[$(i)$] $\xi_1$ is unit spacelike Killing vector field and $\xi_2$, $\xi_3$
are unit timelike Killing vector fields on~$M$.

\item[$(ii)$] $\eta_1$, $\eta_2$, $\eta_3$ are conformal Killing--Yano tensors of
rank~$1$ on $M$.

\item[$(iii)$] $d\eta_1$, $d\eta_2$, $d\eta_3$ are strictly conformal Killing--Yano
tensors of rank~$2$ on $M$.

\item[$(iv)$] $M$ admits Killing--Yano tensors of rank $(2k+1)$, for
$k\in\{0,1,\dots,2n+1\}$.
\end{enumerate}
\end{corollary}
\begin{proof}
$(i)$ From Theorem \ref{4.1} it  follows that $\xi_1$ is a unit
spacelike and $\xi_2$, $\xi_3$ are unit timelike Killing vector f\/ields
on $M$.

$(ii)$ It is well-known that a vector f\/ield is dual to a conformal
Killing--Yano tensor of rank~1 if and only if it is a Killing vector
f\/ield. Consequently, the statement follows from $(i)$, since
$\eta_1$, $\eta_2$, $\eta_3$ are 1-forms dual to the vector f\/ields
$\xi_1$, $\xi_2$, $\xi_3$.

$(iii)$ A direct computation using (\ref{11}), (\ref{12}), (\ref{13})
and (\ref{14}) leads to:
\[
\nabla_Xd\eta_\alpha=-\frac{1}{4n+2}X^*\wedge d^*d\eta_\alpha,\qquad
\alpha\in\{1,2,3\} ,
\]
for any vector f\/ield $X$ on $M$. Therefore, because
$d^2\eta_\alpha=0$, one can write:
\[
\nabla_Xd\eta_\alpha=\frac{1}{3}X\lrcorner
d(d\eta_\alpha)-\frac{1}{4n+2}X^*\wedge d^*(d\eta_\alpha),\qquad
\alpha\in\{1,2,3\} ,
\]
and from (\ref{4}) we deduce that $d\eta_1$, $d\eta_2$, $d\eta_3$ are
conformal Killing--Yano tensors of rank~2 on~$M$. Finally, from
Proposition~\ref{3.4} we deduce that they are strictly conformal
Killing--Yano tensors on~$M$ (i.e.\ they are not Killing tensor
f\/ields); else, if we suppose that $d\eta_\alpha$ is a Killing tensor
f\/ield, then it follows that $\phi_\alpha$ is Killing--Yano, which is
false.

$(iv)$ Since $\eta_1$, $\eta_2$, $\eta_3$ and $d\eta_1$, $d\eta_2$, $d\eta_3$ are
conformal Killing--Yano tensors on $M$, it follows by direct
computations that any linear combination of the forms:
\[
\eta_1\wedge \left(d\eta_1\right)^k, \qquad \eta_2\wedge
\left(d\eta_2\right)^k, \qquad \eta_3\wedge \left(d\eta_3\right)^k  ,
\]
is a Killing--Yano tensor of rank $(2k+1)$, for all
$k\in\{0,1,\dots,2n+1\}$.
\end{proof}

Let us remark that if $M^{4n+3}$ is a manifold endowed with a mixed
3-Sasakian structure
$((\phi_\alpha,\xi_\alpha,\eta_\alpha)_{\alpha=\overline{1,3}},g)$,
then $\xi_1+\xi_2$ and $\xi_1+\xi_3$ are lightlike Killing vector
f\/ields on $M$. Indeed, since $\{\xi_1,\xi_2,\xi_3\}$ are orthogonal
Killing vector f\/ields, with $\xi_1$ unit spacelike vector f\/ield and~$\xi_2$,~$\xi_3$ unit timelike vector f\/ields, it is obvious that
$\xi_1+\xi_2$ and $\xi_1+\xi_3$ are also Killing vector f\/ields and
for $\alpha\in\{2,3\}$ one has:
\begin{gather*}
g(\xi_1+\xi_\alpha,\xi_1+
\xi_\alpha) = g(\xi_1,\xi_1)+2g(\xi_1,\xi_\alpha)+
g(\xi_\alpha,\xi_\alpha)
 = 1+\epsilon_\alpha\nonumber = 0  .
\end{gather*}

Therefore, $\xi_1+\xi_2$ and $\xi_1+\xi_3$ are lightlike Killing
vector f\/ields on $M$.

\begin{corollary}\label{4.3}
Let $(M^{4n+3},(\phi_\alpha,\xi_\alpha,\eta_\alpha)_{\alpha=
\overline{1,3}},g)$ be a manifold endowed with a
mixed $3$-Sasakian structure. Then the distribution spanned by
$\{\xi_1,\xi_2,\xi_3\}$ is integrable and defines a $3$-dimensional
Riemannian foliation on $M$, having totally geodesic leaves of
constant curvature~$1$.
\end{corollary}
\begin{proof}
From (\ref{17}) we obtain:
\[
[\xi_1,\xi_2]=0,\qquad [\xi_2,\xi_3]=2\xi_1,\qquad [\xi_3,\xi_1]=0 ,
\]
and consequently, the distribution spanned by
$\{\xi_1,\xi_2,\xi_3\}$ is integrable and def\/ines a 3-dimen\-sio\-nal
foliation of~$M$. On the other hand, since $\xi_1$, $\xi_2$, $\xi_3$ are
Killing vector f\/ields, it follows that this foliation is Riemannian
(see~\cite{MO}). Moreover, from Theorem \ref{4.1} we conclude that
this foliation has totally geodesic leaves of constant curvature~1.
\end{proof}

\section{Some examples}

We give now examples of manifolds endowed with mixed 3-structures
and mixed 3-Sasakian structures. In particular, we provide some
examples of manifolds which admit Killing--Yano tensors.

\begin{example}
It is easy to see that if we def\/ine
$(\phi_\alpha,\xi_\alpha,\eta_\alpha)_{\alpha=\overline{1,3}}$ in
$\mathbb{R}^3$ by their matrices:
\begin{gather*}
\phi_1=\left(%
\begin{array}{ccc}
  0 & 0 & 1 \\
  0 & 0 & 0 \\
  -1 & 0 & 0
\end{array}%
\right)  ,\qquad
\phi_2=\left(%
\begin{array}{ccc}
  0 & 0 & 0 \\
  0 & 0 & 1 \\
  0 & 1 & 0
\end{array}%
\right)  ,\qquad
\phi_3=\left(%
\begin{array}{ccc}
  0 & -1 & 0 \\
  -1 & 0 & 0 \\
  0 & 0 & 0
\end{array}%
\right) ,\\
\xi_1=\left(%
\begin{array}{c}
  0 \\
  1 \\
  0
\end{array}%
\right)  ,\qquad
\xi_2=\left(%
\begin{array}{c}
  1 \\
  0 \\
  0
\end{array}%
\right)  ,\qquad
\xi_3=\left(%
\begin{array}{c}
  0 \\
  0 \\
  1
\end{array}%
\right) , \\
\eta_1=\left(%
\begin{array}{ccc}
  0 & 1 & 0
\end{array}%
\right)  ,\qquad
\eta_2=\left(%
\begin{array}{ccc}
  -1 & 0 & 0
\end{array}%
\right) ,\qquad
\eta_3=\left(%
\begin{array}{ccc}
  0 & 0 & -1
\end{array}%
\right) ,
\end{gather*}
then $(\phi_\alpha,\xi_\alpha,\eta_\alpha)_{\alpha=\overline{1,3}}$
is a mixed 3-structure on $\mathbb{R}^3$.

We def\/ine now
$(\phi'_\alpha,\xi'_\alpha,\eta'_\alpha)_{\alpha=\overline{1,3}}$ in
$\mathbb{R}^{4n+3}$ by:
\[
\phi'_\alpha=\left(%
\begin{array}{cc}
  \phi_\alpha & 0 \\
  0 & J_\alpha
\end{array}%
\right) , \qquad \xi'_\alpha=\left(%
\begin{array}{c}
  \xi_\alpha \\
  0
\end{array}%
\right)
, \qquad \eta'_\alpha=\left(%
\begin{array}{cc}
  \eta_\alpha & 0 \\
\end{array}%
\right) ,
\]
for $\alpha=1,2,3$, where  $J_1$ is the almost complex structure on~$\mathbb{R}^{4n}$ given by:
\begin{gather}\label{Ianus31}
J_1((x_i)_{i=\overline{1,4n}})=(-x_2,x_1,-x_4,x_3,\dots,-
x_{4n-2},x_{4n-3},-x_{4n},x_{4n-1}) ,
\end{gather}
and $J_2,\ J_3$ are almost product structures on~$\mathbb{R}^{4n}$
def\/ined by:
\begin{gather}\label{Ianus32}
J_2((x_i)_{i=\overline{1,4n}})=
(-x_{4n-1},x_{4n},-x_{4n-3},x_{4n-2},\dots,-x_3,x_4,-x_1,x_2) ,\\
\label{Ianus33}
J_3((x_i)_{i=\overline{1,4n}})=
(x_{4n},x_{4n-1},x_{4n-2},x_{4n-3},\dots,x_4,x_3,x_2,x_1) .
\end{gather}

Since $J_2J_1=-J_1J_2=J_3$, it is easily checked that
$(\phi'_\alpha,\xi'_\alpha,\eta'_\alpha)_{\alpha=\overline{1,3}}$ is
a mixed 3-structure on~$\mathbb{R}^{4n+3}$.
\end{example}

\begin{example}
Let
$(\overline{M},H=(J_{\alpha})_{\alpha=\overline{1,3}},\overline{g})$
be an almost para-hyperhermitian manifold and $(M,g)$ be a
semi-Riemannian hypersurface of $\overline{M}$ with
$g=\overline{g}_{|M}$, having null co-index.

Then, for any $X\in\Gamma(TM)$ and $\alpha\in\{1,2,3\}$, we have the
decomposition: \[J_\alpha X=\phi_\alpha X+F_\alpha X,\] where
$\phi_\alpha X$ and $F_\alpha X$ are the tangent part and the normal
part of $J_\alpha X$, respectively.

But, since $M$ is a semi-Riemannian hypersurface of $M$ with null
co-index, it follows that $TM^\perp=\langle N\rangle $, where $N$ is a unit
space-like vector f\/ield and then $F_\alpha X=\eta_\alpha(X)N,$ where
$\eta_\alpha(X)=\overline{g}(J_\alpha X,N)$.

Consequently, we have the decomposition:
\[
J_\alpha X=\phi_\alpha
X+\eta_\alpha(X)N .
\]

We def\/ine $\xi_\alpha=-J_\alpha N,$ for all $\alpha\in\{1,2,3\}$ and
by straightforward computation it follows that
$((\phi_\alpha,\xi_\alpha,\eta_\alpha)_{\alpha\in\{1,2,3\}},g)$ is a
mixed 3-structure on $M$.

We note that in this way we obtain an entire class of examples of
manifolds which admit mixed 3-structures,  because one can construct
semi-Riemannian hypersurfaces of null co-index, both compact and
non-compact, by using a standard procedure (see \cite{ON}): if
$f\in\mathfrak{F}(\overline{M})$ such that
$\overline{g}({\rm grad}\,f,{\rm grad}\,f)>0$ on $M$ and $c$ is a value of $f$, then
$M=f^{-1}(c)$ is a semi-Riemannian hypersurface of $\overline{M}$.
Moreover, $N={\rm grad}\,f/|{\rm grad}\,f|$ is a unit spacelike vector f\/ield normal
to $M$.

We remark that if $M$ is compact then we can provide examples of
compact hypersurfaces which can be endowed with mixed 3-structures.
In this way we obtain mixed 3-structures on hypersurfaces of the
next compact para-hyper-K\"{a}hler manifolds of dimension 4: complex
tori and primary Kodaira surfaces (see \cite{KMD}).
\end{example}

\begin{example}
Let $(\overline{M},\overline{g})$ be a $(m+2)$-dimensional
semi-Riemannian manifold with index $q\in\{1,2,\dots,m+1\}$ and let
$(M,g)$ be a~hypersurface of $\overline{M}$, with
$g=\overline{g}_{|M}$. We say that $M$ is a~lightlike hypersurface
of $\overline{M}$ if $g$ is of constant rank~$m$ (see~\cite{BD}).
Unlike the classical theory of non-degenerate hypersurfaces, in case
of lightlike hypersurfaces, the induced metric tensor f\/ield~$g$ is
degenerate.

We consider the vector bundle $TM^\perp$ whose f\/ibres are def\/ined
by: \[T_pM^\perp=\{Y_p\in
T_p\overline{M}|\overline{g}_p(X_p,Y_p)=0,\forall \,X_p\in T_pM\},
\qquad \forall \, p\in M.\]

If $S(TM)$ is the complementary distribution of $TM^\perp$ in $TM$,
which is called the screen distribution, then there exists a unique
vector bundle ${\rm ltr}(TM)$ of rank 1 over $M$ so that for any non-zero
section $\xi$ of $TM^\perp$ on a coordinate neighborhood $U \subset
M$, there exists a unique section~$N$ of ${\rm ltr}(TM)$ on $U$
satisfying:
\[
\overline{g}(N,\xi)=1,\qquad
\overline{g}(N,N)=\overline{g}(W,W)=0,\qquad \forall \, W \in
\Gamma(S(TM)_{|U}) .
\]

On another hand, an almost hermitian paraquaternionic manifold is a
triple $(\overline{M},\sigma,\overline{g})$, where $\overline{M}$ is
a smooth manifold, $\sigma$ is a rank 3-subbundle of
${\rm End}(T\overline{M})$ which is locally spanned by an almost
para-hypercomplex structure $H=(J_{\alpha})_{\alpha=\overline{1,3}}$
and $\overline{g}$ is a para-hyperhermitian metric with respect to~$H$.

In \cite{IMV} it is proved that there is a mixed 3-structure on any
lightlike hypersurface $M$ of an almost hermitian paraquaternionic
manifold $(\overline{M},\sigma,\overline{g})$, such that $\xi$ and
$N$ are globally def\/ined on~$M$. In particular, the above conditions
are satisf\/ied if we consider the lightlike hypersurface~$M$ of
$(R_{4}^{8},\overline{g},\{J_1,J_2,J_3\})$ def\/ined by (see
\cite{IMV2}):
\[
f((t_i)_{i=\overline{1,7}})=
(t_1,t_2,t_1+t_3,t_4+t_5,t_3+t_5+t_6,t_6+t_7,t_7,t_2) ,
\]
where the structures $J_1$, $J_2$, $J_3$ and the metric
$\overline{g}$ on $R_{4}^{8}$ are given by:
\begin{gather*}
\overline{g}((x_i)_{i=\overline{1,8}},(y_i)_{i=\overline{1,8}})=
-\sum_{i=1}^{4}x_iy_i+\sum_{i=5}^{8}x_iy_i ,\\
J_1((x_i)_{i=\overline{1,8}})=(-x_2,x_1,-x_4,x_3,-x_6,x_5,-x_8,x_7),\\
J_2((x_i)_{i=\overline{1,8}})=(-x_7,x_8,-x_5,x_6,-x_3,x_4,-x_1,x_2),\\
J_3((x_i)_{i=\overline{1,8}})=(x_8,x_7,x_6,x_5,x_4,x_3,x_2,x_1).
\end{gather*}
\end{example}

\begin{example}
The unit pseudo-sphere $S^{4n+3}_{2n+1}$ is the canonical example of
manifold with a mixed 3-Sasakian structure \cite{IMV3}. This
structure is obtained by taking $S^{4n+3}_{2n+1}$ as hypersurface of
$(\mathbb{R}^{4n+4}_{2n+2},\overline{g})$. It is easy to see that on
the tangent spaces $T_pS^{4n+3}_{2n+1}$, $p\in S^{4n+3}_{2n+1}$, the
induced metric $g$ is of signature $(2n+1,2n+2)$.

If $(J_\alpha)_{\alpha=\overline{1,3}}$ is the canonical
para-hypercomplex structure on the $\mathbb{R}^{4n+4}_{2n+2}$ given
by (\ref{Ianus31})--(\ref{Ianus33}) and $N$ is the unit spacelike
normal vector f\/ield to the pseudo-sphere, we can def\/ine three vector
f\/ields on $S^{4n+3}_{2n+1}$ by: \[\xi_\alpha=-J_\alpha N,\qquad
\alpha=1,2,3.\]

If $X$ is a tangent vector to the pseudo-sphere then $J_\alpha X$
uniquely decomposes onto the part tangent to the pseudo-sphere and
the part parallel to $N$. Denote this decomposition by:
\[
J_\alpha X=\phi_\alpha X+\eta_\alpha(X)N .
\]

This def\/ines the 1-forms $\eta_\alpha$ and the tensor f\/ields
$\phi_\alpha$ on $S^{4n+3}_{2n+1}$, where $\alpha=1,2,3$. Now we can
easily see that
$(\phi_\alpha,\xi_\alpha,\eta_\alpha)_{\alpha=\overline{1,3}}$ is a
mixed 3-Sasakian structure on $S^{4n+3}_{2n+1}$.
\end{example}

\begin{example}
Since we can recognize the unit pseudo-sphere $S^{4n+3}_{2n+1}$ as
the projective space $P^{4n+3}_{2n+1}(\mathbb{R})$, by identifying
antipodal points, we also have that $P^{4n+3}_{2n+1}(\mathbb{R})$
admits a mixed 3-Sasakian structure.
\end{example}

Applying Corollary \ref{4.2} we deduce that both $S^{4n+3}_{2n+1}$
and $P^{4n+3}_{2n+1}(\mathbb{R})$ admit spacelike, timelike and
lightlike Killing vector f\/ields, conformal Killing--Yano tensors of
rank~1, strictly conformal Killing--Yano tensors of rank 2 and
Killing--Yano tensors of rank $(2k+1)$, for $k\in\{0,1,\dots,2n+1\}$.

\subsection*{Acknowledgements}

We would like to thank the referees for carefully reading the paper
and making valuable comments and suggestions. This work was
supported by CNCSIS Programs, Romania.

\pdfbookmark[1]{References}{ref}
\LastPageEnding

\end{document}